\documentclass[11pt]{article}

\usepackage[margin=1.08in]{geometry}
\usepackage{amsmath,amssymb,amsthm,mathtools}
\usepackage{booktabs}
\usepackage{graphicx}
\usepackage[expansion=false]{microtype}
\usepackage{xcolor}
\usepackage{hyperref}
\usepackage{doi}
\usepackage[nameinlink,capitalize,noabbrev]{cleveref}

\hypersetup{
  colorlinks=true,
  linkcolor=blue!55!black,
  citecolor=blue!55!black,
  urlcolor=blue!60!black,
  pdftitle={$144$ circles tangent to three conics},
  pdfauthor={Taylor Brysiewicz}
}

\newtheorem{theorem}{Theorem}[section]

\usepackage{listings}
\definecolor{maroon}{RGB}{133, 5, 63}
\definecolor{teal}{RGB}{0, 128, 96}
\definecolor{forestgreen}{RGB}{34, 139, 34}
\lstdefinelanguage{julia}{
basicstyle=\small\ttfamily,
alsoletter=",
classoffset=1,
comment=[l]{\#},
commentstyle=\color{gray},
keywords={solve, differentiate, subs, sub, real_solutions, det, sum, max, conditional_count, map, filter, map_filter, union, zip, product, is_finite, is_successful, is_nonsingular, track, is_real, is_success, first, evaluate, flatten, bitmask_filter, accumulate, polynomial_interpolants, reverse, coefficients, vcat, solution_from_necklace, length, prod, compress, System, certify total_degree_start_solutions, degrees, variables, iterate, struct, collect, convert, stretched_cube, stretched_cubes, weight_vector, weight_vectors, perm_to_segments, perm_to_mixedcell, mixed_cell_iterator, rand_approx_unit, norm, fixed, polyhedral_system, reduce, findall, println, eachrow, eval, parse, readlines, certify, certificates, real},
keywordstyle={\color{teal}},
classoffset=2,
morekeywords={@var, @time, for, end, if, while, else, begin, assert, global, in, all, any },
keywordstyle={\color{maroon}},
classoffset=3,
morekeywords={using, function, return, const},
keywordstyle={\color{blue}},
classoffset=4,
morekeywords={julia, >},
keywordstyle={\color{forestgreen}},
xleftmargin=0.2cm,
xrightmargin=1em,
columns=fullflexible,
keepspaces=true,
}

\title{\textbf{144 real circles tangent to three conics}}

\author{Taylor Brysiewicz}

\date{}

\begin{document}

\maketitle

\begin{abstract}
\hspace{-23pt} \noindent We exhibit three conics with $144$ tritangent circles, disproving the conjectured maximum of~$136$. 
\end{abstract}

\section{Statement of Main Result}

We disprove \cite[Conjecture 1.4]{BreidingLindberg} by counterexample and numerical certification (see \hyperref[appendix]{Appendix}):

\begin{theorem}
\label{thm:maintheorem}
There exists a triple of conics admitting $144$ tritangent circles.
\end{theorem}

 \begin{figure}[!htpb]
 \includegraphics[scale=0.155]{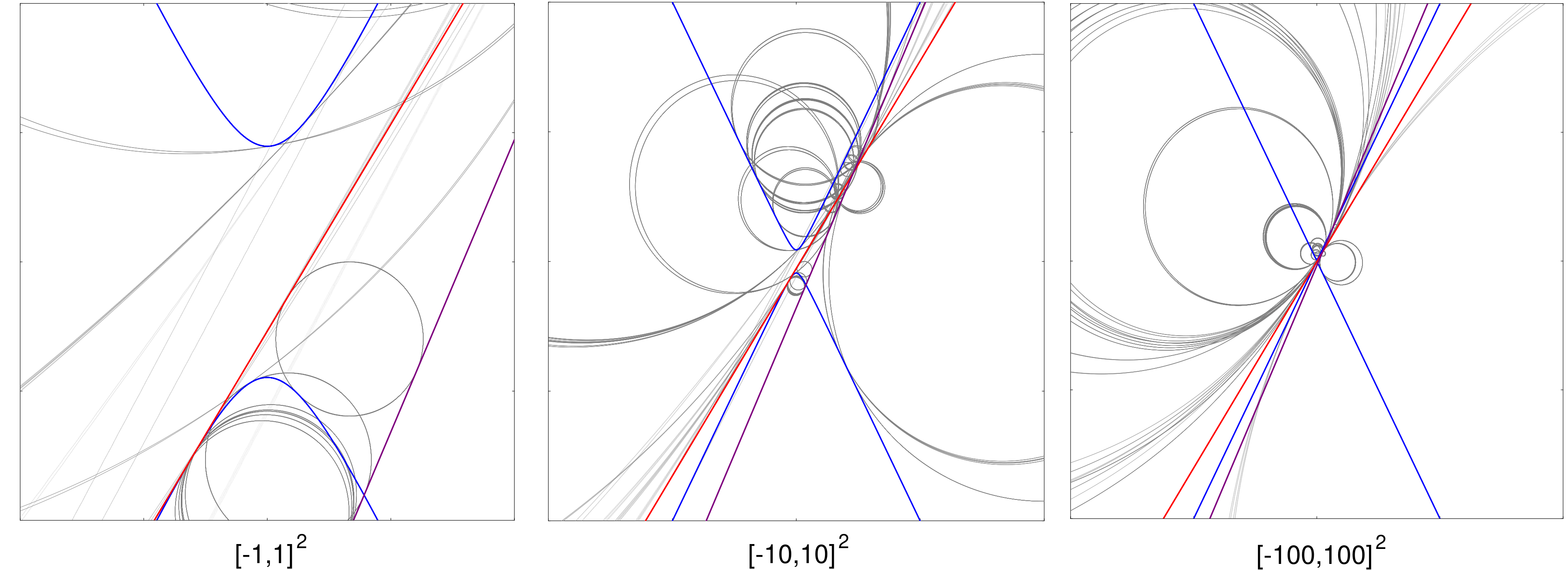} \vspace{-5pt}
 \caption{Three windows displaying three conics and $144$ tritangent circles to them.}
 \end{figure}

\noindent The three conics $\textbf{C}=(C_1,C_2,C_3)$ which witness this result are the respective zero sets of  
\begin{equation}
\label{eq:conics}
\hspace{-35pt}\textbf{Q} = {\small{
\begin{cases}
Q_1(x,y)
=& 20x^2-5y^2+1, \\[3pt]
Q_2(x,y)
=&
\left(
y-\frac{2399}{1000}
-\frac{16009}{10000}\left(x-\frac{833}{500}\right)
\right)
\left(
y-\frac{2399}{1000}
-\frac{16011}{10000}\left(x-\frac{833}{500}\right)
\right)
+10^{-12}, \text{ and } \\[3pt]
Q_3(x,y)
=&
\left(
y-\frac{3857}{1000}
-\frac{22609}{10000}\left(x-\frac{2497}{1000}\right)
\right)
\left(
y-\frac{3857}{1000}
-\frac{22611}{10000}\left(x-\frac{2497}{1000}\right)
\right)
+10^{-12}. 
\end{cases}}}\end{equation}

\noindent The conics $C_2$ and $C_3$ are small perturbations of  
double lines:
\[
Q_2
=
L_1^2-\delta^2(x-x_1)^2+\varepsilon,
\qquad
Q_3
=
L_2^2-\delta^2(x-x_2)^2+\varepsilon,
\]
where $L_i=y-y_i-m_i(x-x_i)$ for $i=1,2$ with
\[
(x_1,y_1)
= {\footnotesize{
\left(\frac{833}{500},\frac{2399}{1000}\right)}},
\qquad
(x_2,y_2)
={\footnotesize{
\left(\frac{2497}{1000},\frac{3857}{1000}\right)}},
\]
and
\[
m_1={\footnotesize{\frac{1601}{1000}}},
\qquad
m_2={\footnotesize{\frac{2261}{1000}}},
\qquad
\delta=10^{-4},
\qquad
\varepsilon=10^{-12}.
\]
Thus, as $\delta,\varepsilon\to0$, the conics $C_2$ and $C_3$
degenerate to the double lines $L_1^2$ and $L_2^2$, respectively.
\newpage

For three general complex conics there are $184$ tritangent circles
\cite{BreidingLindberg}. Determining the maximum number of real solutions
in a parametrized polynomial family is generally difficult. Some computational
approaches move through parameter space in search of chambers with more real
solutions \cite{Dietmaier, GriffinHauenstein}. A complementary strategy in
real enumerative geometry is to degenerate the constraints until the problem
decomposes into simpler enumerative problems, make those constituent problems
simultaneously as real as possible, and perturb to a nearby instance
\cite{BreidingSturmfelsTimme,TangentQuadrics,RTV,SottileEffective}.

The difficulty in successfully applying a degeneration technique is one of
\textit{balance}. Choosing a family of parameters that is too special may
nontrivially restrict the number of real solutions. Choosing
a family that is not special enough may fail to offer a sufficiently fine
decomposition for which the constituent parts can be understood. We improve
upon the previous record of $136$ real tritangent circles obtained in
\cite{BreidingLindberg,thesis} exactly by relaxing the triangle degeneration
of \Cref{secsec:triangledecomposition} to one which keeps one of the three
conics irreducible. This degeneration was explored in \cite{thesis}.

 \textbf{AI's role in the discovery:} The thesis \cite{thesis} was written in French, not mentioned in \cite{BreidingLindberg}, and not known to the present author. OpenAI's GPT-5.6 Sol model found this source and identified its prescribed degeneration as relevant to the search for more than $136$ real solutions. Guided by this degeneration, the author searched numerically for compatible real instances of the four constituent enumerative problems, leading to the configuration above. The author used the \texttt{julia} package \texttt{HomotopyContinuation.jl} \cite{HC} to solve the associated system, however, several solutions were so numerically ill-conditioned that the software reported them as singular. Nonetheless, numerical certification \cite{certify} in that same package certified that all of the $184$ solutions to this system are isolated and distinct, that $144$ are real and represent $144$ distinct circles.

 \section*{Acknowledgements}
TB was supported by NSERC Discovery Grant RGPIN-2023-03551. 

 \section{Two Degenerations}
 We describe two degenerations. The first degenerates the three conics to double lines. The second is a relaxation which only degenerates two of the three conics to double lines. 
 \subsection{The triangle degeneration}
 \label{secsec:triangledecomposition}
 The pentagon degeneration of Fulton (see \cite{BreidingSturmfelsTimme}), and independently Ronga, Tognoli, and Vust \cite{RTV}, degenerates five conics to a pentagon of double lines, thus degenerating the  $3264$ conics tangent to them. They form groups of $32$ in the limit, which are further partitioned into solutions of simpler enumerative problems regarding conics tangent to fixed lines and passing through fixed points. These simple enumerative problems can be totally real, and moreover, can be simultaneously made real by a set of five lines  with fixed marked points. Once one finds a simultaneously totally real pentagon, the idea is to search nearby to find nondegenerate totally real instances of the original enumerative problem. Fulton first observed that all $3264$  solutions could be real via this degeneration and communicated the argument to Sottile \cite[Chapter~9]{SottileBook}.  Independently, Ronga–Tognoli–Vust gave a detailed published proof in $1997$ \cite{RTV}. Much later, Breiding, Sturmfels, and Timme
\cite{BreidingSturmfelsTimme} produced an explicit numerical instance of five
conics with all $3264$ tangent conics real and rigorously certified the solutions.

\begin{figure}[!htpb]
\begin{center}
\includegraphics[scale=0.2]{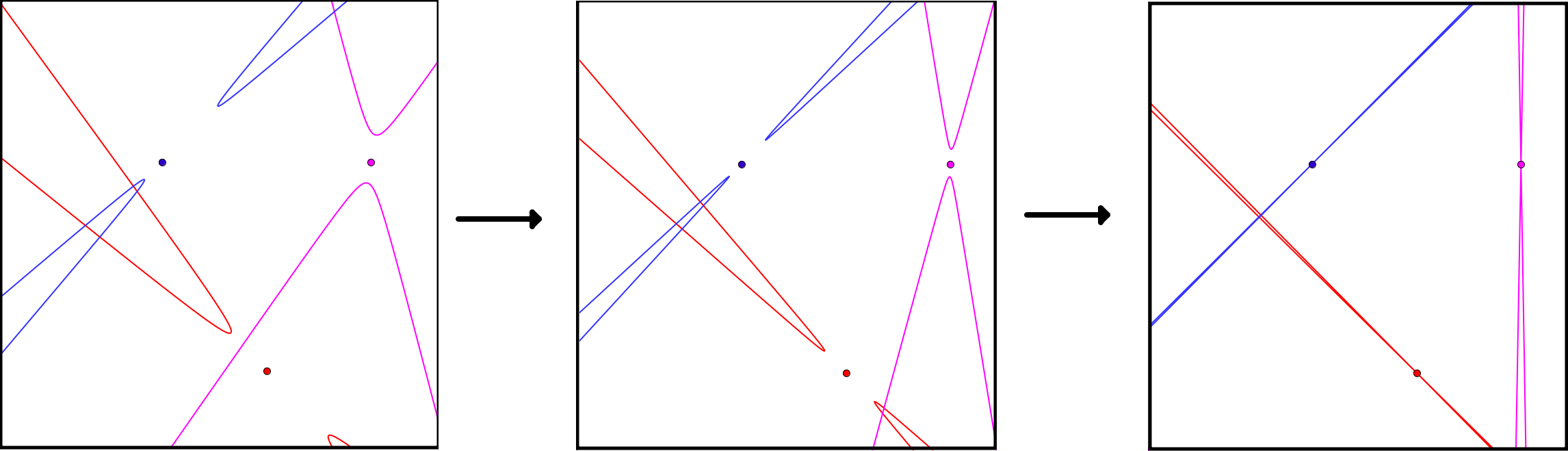}
\end{center}
\caption{A visualization of the triangle degeneration. The hyperbolas on the left degenerate to double lines on the right,
with marked points indicating the limiting contacts of their branches.}
\label{fig:triangle}
\end{figure}

 The triangle analogue of this technique degenerates three conics to a triangle of double lines (see \Cref{fig:triangle}) making the $184$ circles tangent to them collide in $23$ groups of $8$, each satisfying a smaller enumerative problem. Specifically, these are the problems $\{P^{3-k}L^{k}\}_{k=0}^3$ of finding circles tangent to $k=0,\ldots,3$  lines and through $3-k$ points. The counts to these problems are $1,2,4,$ and $4$, and there are ${{3}\choose{k}}$ of each, respectively. This degeneration thus interprets $184$ as
 \[
 184 = 8\left({{3}\choose{0}}\#P^3+ {{3}\choose{1}} \cdot \#P^2L + {{3}\choose{2}} \cdot \#PL^2 + {{3}\choose{3}}\cdot \#L^3\right)=8\left(1 + 3 \cdot 2 + 3 \cdot 4 +  4\right).
 \]
 
Vilmart proved in \cite[Thm.~1.1.1]{thesis} that this degeneration
produces a nearby triple of smooth hyperbolas with $136$ real
tritangent circles.  In the corollary to the proof of
\cite[Thm.~1.1.1]{thesis}, she further observes that $136$ is the
largest real count obtainable by this particular degeneration.
The same degeneration was rediscovered in \cite{BreidingLindberg}
thirteen years later.  Since at most  two of the four complex circles through a point and tangent to two lines can be real, any triple of conics nearby a triangle of double lines can have at most $8\cdot(1+3\cdot 2 + 3 \cdot 2 + 4) = 136$ real tritangent circles. 

 \subsection{A mixed degeneration} In the ``Perspectives'' of her thesis \cite[pp.~131--132]{thesis},
Vilmart proposes a less degenerate construction in which one hyperbola
is kept while the other conics degenerate to double lines. As with the triangle degeneration, the enumerative problem of degree $184$ decomposes, but this time as \begin{equation}
\label{eq:thesisdegeneration}
184
=
4\left(\#QP^2+\#QPL+\#QLP+\#QL^2\right)
=
4(6+12+12+16).
\end{equation}

For the limiting configuration underlying our choice of conics \eqref{eq:conics}, the four constituent problems have  $
 6, 8, 8,$ and $14$
 real solutions, respectively. Consequently, the corresponding real count is $
 144 = 4(6+8+8+14)$. We remark that this was obtained despite the degree $12$ problems not being as real as possible.

 \begin{table}[!htpb]
\centering {\footnotesize{
\begin{tabular}{l|c|ccc}
\toprule
& $\mathbb{C}$ & & $\mathbb{R}$ &\\ \hline 
Subproblem
& \shortstack{\\Count}
& \shortstack{\\ Individually \\ Achieved}
& \shortstack{Upper Bound}
& \shortstack{Simultaneously\\  achieved} \\
\midrule
$QP^2$
& $6$
& $6$
& $6$
& $6$ \\

$QPL$
& $12$
& $10$
& $10$
& $8$ \\

$QLP$
& $12$
& $10$
& $10$
& $8$ \\

$QL^2$
& $16$
& $14$
& $16$
& $14$ \\
\midrule
Sum
& $46$
& $40$
& $42$
& $36$ \\

$ \times 4$
& $184$
& $160$
& $168$
& $144$ \\
\bottomrule
\end{tabular}}}
\caption{Complex and real solution counts for the four constituent problems in the degeneration.}
\label{tab:numerology}
\end{table} 
The relevant numerology is summarized in \autoref{tab:numerology}.
Vilmart proves in \cite[Thm.~4.0.2]{thesis} that $QPL$ has at most ten
real solutions, and gives a construction realizing all ten
in \cite[Sec.~5.2]{thesis}.  She considers the mixed degeneration
above in \cite[pp.~131--132]{thesis}.  Assuming that her candidate counts $6$, $10$, $10$, and $12$ are simultaneously achievable, she estimates a possible real count of 
\[
152=4(6+10+10+12).
\]
She does not claim this as a realized configuration.  She does,
however, observe that the construction of Chapter~5 can be arranged
so that the first three constituent problems are as real as possible, contributing $4 \cdot (6 + 10 + 10) = 104$ to the total real count, with multiplicity.  What remains unresolved is the $QL^2$ constituent and
whether the degenerate configuration can be perturbed while maintaining this real count.

For $QL^2$, Vilmart explicitly poses
the maximal real count as an open problem
\cite[p.~132]{thesis}. A real count of $12$ was known to be possible, giving an estimate of $152$ possible real solutions. Our new record of $14$
therefore improves the previously known count, making $160=4(6+10+10+14)$
a natural next target for the problem of degree $184$.

\section{Visualizing the $144$-chamber}
We use \texttt{AdaptiveVisualization.jl}  \cite{AdaptiveVisualizationJL} to visualize the parameter space of triples of conics near $\textbf{Q}=(Q_1,Q_2,Q_3)$ as in  \eqref{eq:conics}. Specifically, we visualize the $(t_1,t_2)$-plane $\Lambda:$
\begin{equation}
\label{eq:parametrization}
\Lambda: \textbf{Q} + 10^{-12} t_1 \textbf{P}_1 + 10^{-12} t_2 \textbf{P}_2
\end{equation}
containing $\textbf{Q}$ and parametrized by
{\footnotesize{
\begin{align*}
\hspace{-6pt}\textbf{P}_1 \hspace{-3pt}&=\hspace{-3pt} \begin{pmatrix}
-0.2680415043319x\hspace{-2pt}-\hspace{-2pt}0.9125539464648x^2
\hspace{-2pt}-\hspace{-2pt}0.0604015811180y\hspace{-2pt}-\hspace{-2pt}1.2431768299194y^2
\hspace{-2pt}-\hspace{-2pt}1.1954903428489xy
\\
-0.0503305778780x\hspace{-2pt}-\hspace{-2pt}0.4561730742749x^2
\hspace{-2pt}-\hspace{-2pt}0.7755220919526y\hspace{-2pt}-\hspace{-2pt}0.1692096889321y^2
\hspace{-2pt}-\hspace{-2pt}1.5692576120810xy
\\
-0.6819277495837x\hspace{-2pt}+\hspace{-2pt}0.2921430259032x^2
\hspace{-2pt}-\hspace{-2pt}0.3974300425398y\hspace{-2pt}-\hspace{-2pt}0.1570871803620y^2
\hspace{-2pt}+\hspace{-2pt}1.0685570773946xy
\end{pmatrix} \\ 
\textbf{P}_2 &= \begin{pmatrix}
0.1347246634245x\hspace{-2pt}-\hspace{-2pt}0.7829058032091x^2
\hspace{-2pt}-\hspace{-2pt}1.8873834159597y\hspace{-2pt}-\hspace{-2pt}0.0532721015859y^2
\hspace{-2pt}+\hspace{-2pt}0.7698217491853xy
\\
-0.5023516834741x\hspace{-2pt}+\hspace{-2pt}1.2872693960393x^2
\hspace{-2pt}+\hspace{-2pt}0.5395746942042y\hspace{-2pt}-\hspace{-2pt}2.0600574420752y^2
\hspace{-2pt}-\hspace{-2pt}0.7384063387941xy
\\
-0.5254486877266x\hspace{-2pt}+\hspace{-2pt}0.6581332367453x^2
\hspace{-2pt}+\hspace{-2pt}0.2066888399141y\hspace{-2pt}+\hspace{-2pt}0.1663650392275y^2
\hspace{-2pt}-\hspace{-2pt}1.1195775945642xy
\end{pmatrix}
\end{align*}}}

\begin{figure}[!htpb]
\begin{center}
\includegraphics[scale=0.8]{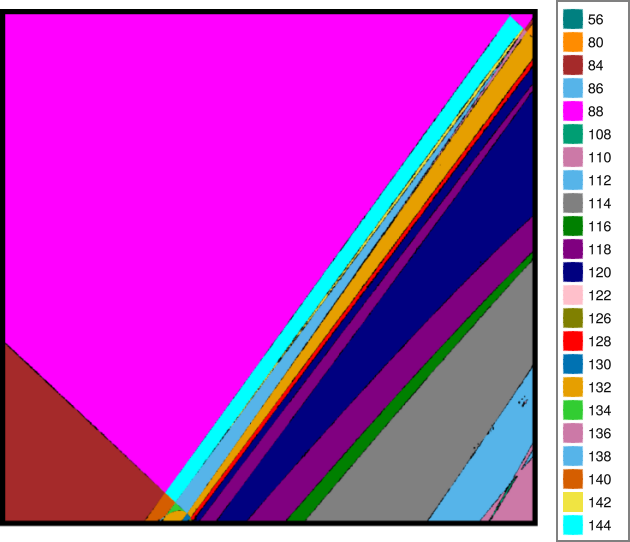}
\end{center}
\caption{The $(t_1,t_2)$-parameter space in window $[-188,9.8] \times[-187,10.25]$. This view contains the entire $\Lambda$-slice of the $144$ chamber.}
\label{fig:fullchamber}
\end{figure}

\begin{figure}[!htpb]
\begin{center}
\includegraphics[scale=0.65]{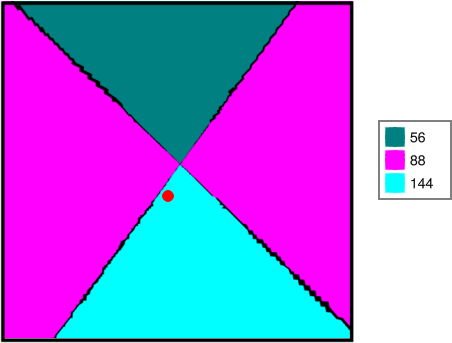}
\end{center}\caption{Near the top of the $144$ chamber, with
$(0,0)$ corresponding to $\mathbf Q$ marked in red.}
\label{fig:top}
\end{figure}

\begin{figure}[!htpb]
\begin{center}
\includegraphics[scale=0.8]{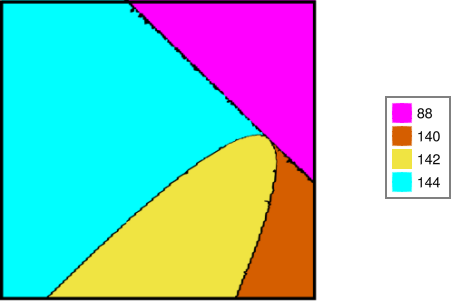}
\end{center}
\caption{Successively zooming in on the right corner of the $144$ chamber.}
\label{fig:right}
\end{figure}

\begin{figure}[!htpb]
\begin{center}
\includegraphics[scale=1.1]{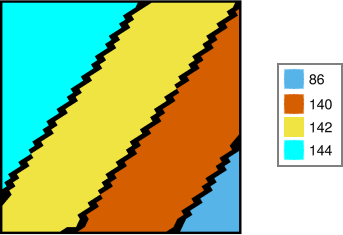}
\end{center}
\caption{Successively zooming in on the right edge of the $144$ chamber on the scale of $0.02$.}
\label{fig:rightedge}
\end{figure}

\begin{figure}[!htpb]
\begin{center}
\includegraphics[scale=0.6]{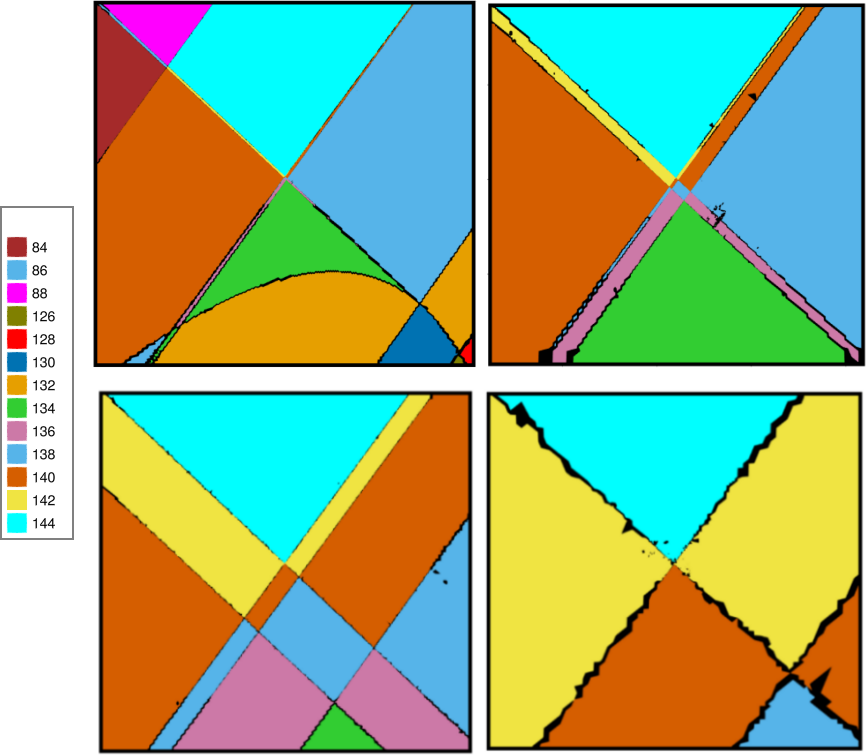}
\end{center}
\caption{Successively zooming in on the bottom corner $(-123.395,-179.491)$.}
\label{fig:bottom}
\end{figure}

\begin{figure}[!htpb]
\begin{center}
\includegraphics[scale=0.7]{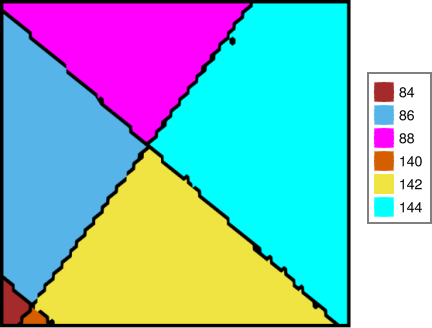}
\end{center}
\caption{The left corner of the $144$ chamber.}
\label{fig:left}
\end{figure}

The remainder of this paper is a gallery of images of the chamber containing $\textbf{Q}$ within $\Lambda$. That chamber is colored cyan. It seems to be roughly diamond-shaped. We showcase 
\begin{itemize}
\item the full chamber \autoref{fig:fullchamber}, 
\item the top corner \autoref{fig:top} near the degenerate triple, 
\item the right corner \autoref{fig:right} where there seems to be a parabolic boundary, 
\item the right edge \autoref{fig:rightedge}, where zooming-in reveals thin adjacent chambers,
\item the bottom corner \autoref{fig:bottom}, again, with very small adjacent chambers, and 
\item the left corner \autoref{fig:left} which also reveals small adjacent chambers. 
\end{itemize}

\begin{figure}[!htpb]
\begin{center}
\includegraphics[scale=0.6]{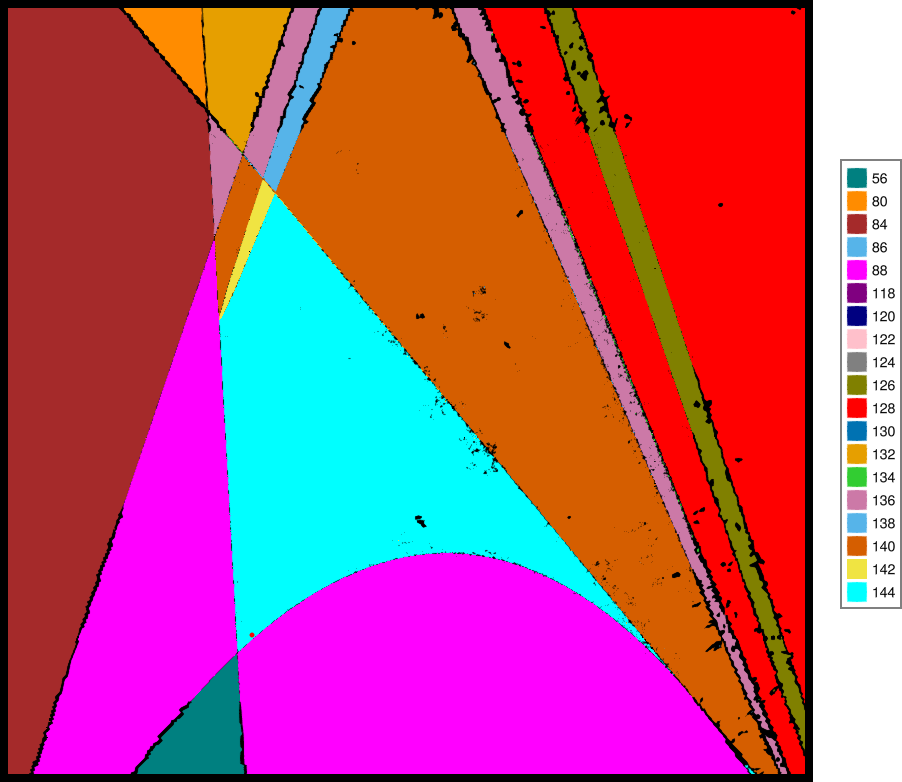}
\end{center}
\caption{A second two-dimensional slice $\Lambda'$ through a
neighborhood of $\mathbf Q$. The red dot marks $\mathbf Q$, near the
apparent boundary between regions with $144$ and $88$ real solutions. Black denotes unresolved evaluations. }
\label{fig:anotherview}
\end{figure}

Finally, we give another $2$-plane intersection illustration. This view is shown in \Cref{fig:anotherview}, although our software had more difficulty reliably computing the $184$ solutions within this slice (indicated by many small black regions).

\newpage

 \section*{Appendix}
\label{appendix} The code below numerically certifies (see \cite{certify}) that the
corresponding marked tangency system has $184$ isolated and distinct
solutions. Of these, $144$ are real, all have positive squared radius,
and their $(s,t,r)$-coordinates are pairwise distinct. Thus they
determine $144$ distinct real circles, proving \Cref{thm:maintheorem}. Note that our formulation fixes the constant terms of the input conics to be $1$, explaining the apparent difference to \eqref{eq:conics}.
\subsection*{Code}
\begin{lstlisting}[language=julia]
monodromy_seed = 0x00000144; tracking_seed  = 0x14400184
using HomotopyContinuation, LinearAlgebra, Arblib
#Construct parametrized system
@var x y s t r u[1:3] v[1:3] a[1:3,1:5]
Q = [a[i,1]*x^2 + a[i,2]*x*y + a[i,3]*y^2 + a[i,4]*x + a[i,5]*y + 1 for i in 1:3]
circle = (x-s)^2 + (y-t)^2 - r
f = [evaluate(circle,[x,y]=>[u[i],v[i]]) for i in 1:3]
g = [evaluate(Q[i],[x,y]=>[u[i],v[i]]) for i in 1:3]
h = [det(hcat(differentiate(f[i],[u[i],v[i]]),
              differentiate(g[i],[u[i],v[i]]))) for i in 1:3]
F = System(vcat(f,g,h),
    variables=[u[1],v[1],u[2],v[2],u[3],v[3],s,t,r],
    parameters=vec(a))
# Exact target parameters
p144 = Rational{BigInt}[
    20, 1780000687500//49976818751, 170404033000000//106650281457997,
    0, -20012500000000//449791368759, -452200000000000//319950844373991,
    -5, 6250000000000//449791368759, 100000000000000//319950844373991,
    0, -5368673116750//449791368759, -808857822406000//319950844373991,
    0, 1117775000000//149930456253, 119247800000000//106650281457997
]
# Compute one generic fiber of degree 184 by monodromy
M = monodromy_solve(F; seed = monodromy_seed,  target_solutions_count=184)
# Track those 184 solutions to the exact target parameter
R = solve(F,solutions(M); start_parameters=parameters(M), 
          target_parameters=p144, seed = tracking_seed)
PR = path_results(R)
S = [p.solution for p in PR];
#The following line certifies 144/184 real/complex solutions
C = certify(F,S,p144; max_precision=8192)
DC = distinct_certificates(C)
@assert length(DC) == 184
RC = filter(is_real, DC)
@assert length(RC) == 144
@assert all(c -> Arblib.is_positive(
    real(certified_solution_interval(c)[9])), RC)
# Certify that the 144 incidence solutions give distinct circles
X = certified_solution_interval.(RC)
disjoint(a,b) = Arblib.is_positive(a-b) || Arblib.is_positive(b-a)
@assert all(
    any(disjoint(real(X[i][k]),real(X[j][k])) for k in 7:9)
    for i in 1:143 for j in i+1:144)
\end{lstlisting}

\newpage

\bibliographystyle{plain}
\bibliography{mybib}

\noindent
\textsc{Taylor Brysiewicz}\\
Department of Mathematics, Western University\\
London, Ontario N6A 5B7, Canada\\
\textit{Email address:} \texttt{tbrysiew@uwo.ca}
\end{document}